\documentclass[preprint,1p]{elsarticle}

\makeatletter
 \def\ps@pprintTitle{%
 	\let\@oddhead\@empty
 	\let\@evenhead\@empty
 	\def\@oddfoot{\footnotesize\itshape
 		{} \hfill\today}%
 	\let\@evenfoot\@oddfoot
 }
\makeatother
\usepackage[unicode]{hyperref}
\usepackage[makeroom]{cancel}
 
\usepackage{soul}
\usepackage{tikz}

\usepackage{latexsym}
\usepackage{indentfirst}
\usepackage{amsxtra}
\usepackage{amssymb}
\usepackage{amsthm}
\usepackage{amsmath}
\usepackage{mathrsfs} 

\usepackage{xcolor}
\usepackage{color}

\usepackage{amsfonts}

\usepackage[capitalise]{cleveref}

\usepackage[capitalise]{cleveref}

\newtheorem{theor}{Theorem}
\newtheorem*{theor*}{Theorem}
\newtheorem{prop}[theor]{Proposition}
\newtheorem{lemma}[theor]{Lemma}
\newtheorem{cor}[theor]{Corollary}
\newtheorem*{cor*}{Corollary}
\theoremstyle{definition}               
\newtheorem{defin}[theor]{Definition}
\newtheorem{ex}{Example}

\newtheorem{rem}{Remark}
\newtheorem{rems}[rem]{Remarks}
\newtheorem{que}{Question}

\DeclareMathOperator{\Hom}{Hom}
\DeclareMathOperator{\Sym}{Sym}

\DeclareMathOperator{\End}{End}
\DeclareMathOperator{\id}{id}
\DeclareMathOperator{\im}{Im}

\DeclareMathOperator{\Map}{Map}
\DeclareMathOperator{\Soc}{Soc}

\DeclareMathOperator{\E}{E}

\newcommand{\s}[1]{S_{#1}}
\newcommand{\phii}[2]{\phi_{#1,#2}}

\begin{document}

\begin{frontmatter}

\title{Rota--Baxter operators on Clifford semigroups and the Yang--Baxter equation}
	\tnotetext[mytitlenote]{This work was partially supported by the Dipartimento di Matematica e Fisica ``Ennio De Giorgi'' - Università del Salento. The authors are members of GNSAGA (INdAM).
	}

\author{Francesco CATINO}
\ead{francesco.catino@unisalento.it}

\author{Marzia MAZZOTTA
}
\ead{marzia.mazzotta@unisalento.it}

\author{Paola STEFANELLI
}
\ead{paola.stefanelli@unisalento.it}


\address{Dipartimento di Matematica e Fisica ``Ennio De Giorgi", Universit\`a del Salento, Via Provinciale Lecce-Arnesano, 73100 Lecce (Italy)}

\begin{abstract}
In this paper, we introduce the theory of Rota--Baxter operators on Clifford semigroups, useful tools for obtaining dual weak braces, i.e., triples $\left(S,+,\circ\right)$ where $\left(S,+\right)$ and $\left(S,\circ\right)$ are Clifford semigroups such that
$a\circ\left(b+c\right) = a\circ b - a +a\circ c$ and $a\circ a^- = -a+a$, for all $a,b,c\in S$. To each algebraic structure is associated a set-theoretic solution of the Yang--Baxter equation that has a behaviour near to the bijectivity and non-degeneracy. Drawing from the theory of Clifford semigroups, we provide methods for constructing dual weak braces and deepen some structural aspects, including the notion of ideal.
\end{abstract}

\begin{keyword}
 Rota--Baxter operator \sep Quantum Yang--Baxter equation \sep Set-theoretic solution \sep Inverse semigroup\sep Clifford semigroup \sep Skew brace  \sep Brace \sep Weak brace
\MSC[2020] 17B38 \sep 16T25\sep 81R50\sep 20M18
\end{keyword}
\end{frontmatter}

 \section*{Introduction}

Rota--Baxter operators on commutative algebras firstly appeared in 1960 in G. Baxter's probability studies \cite{Bax60} and were subsequently explored by several authors, including Atkinson \cite{At63}, Cartier \cite{Cartier72}, and Rota \cite{Ro69, Ro95}, among others. An entire monograph has been dedicated by Guo \cite{Gu12} on this topic. Moreover, Guo, Lang, and Sheng \cite{GuLaSh20} have introduced the notion of Rota--Baxter operator on groups, pursued further by Bardakov and Gubarev \cite{BaGu21groupsx}. Recently, the latter have found in \cite{BaGu21} a strict connection between Rota--Baxter operators on groups and the Yang--Baxter equation, a fundamental identity of theoretical physics that takes the name from Yang \cite{Ya67} and R. Baxter \cite{Ba72}, and of which Drinfel'd \cite{Dr92} posed the problem of finding all set-theoretic solutions. Specifically, if $S$ is a non-empty set, a map
$r:S\times S\to S\times S$ is a \emph{set-theoretic solution of the Yang--Baxter equation}, or shortly a \emph{solution}, if it  satisfies the relation
\begin{align*}
\left(r\times\id_S\right)
\left(\id_S\times r\right)
\left(r\times\id_S\right)
= 
\left(\id_S\times r\right)
\left(r\times\id_S\right)
\left(\id_S\times r\right).
\end{align*}
Writing a solution $r$ as $r\left(x,y\right) = \left(\lambda_x\left(y\right), \rho_y\left(x\right)\right)$, with $\lambda_x, \rho_x$ maps from $S$ into itself,  then $r$ is \emph{left non-degenerate} if $\lambda_x\in \Sym_S$, \emph{right non-degenerate} if $\rho_x\in \Sym_S$, for every $x\in S$,  \emph{non-degenerate} if it is both left and right non-degenerate. 

As proved in \cite[Proposition 3.1]{BaGu21}, every Rota--Baxter operator on a group determines a skew brace, an algebraic structure introduced by Guarnieri and Vendramin in 2017 that always gives rise to a bijective non-degenerate solution (see \cite[Theorem 3.1]{GuVe17}).  

In this regard, Caranti and Stefanello \cite{CaSt22x} have recently described gamma functions that characterise skew braces obtained from Rota-Baxter operators on groups.

A generalisation of the skew brace that determines a solution has been newly introduced in \cite{CaMaMiSt21x}, and it is named weak brace. To give the definition, we firstly recall that a \emph{Clifford semigroup} $(S, \cdot)$ is an inverse semigroup, i.e., a semigroup for which every $a \in S$ admits a unique $a^{-1} \in S$ satisfying $aa^{-1} a=a$ and $a^{-1} a^{-1}=a^{-1}$, having central idempotents. 
Into detail, a triple $(S, +, \circ)$ is a \emph{weak (left) brace} if $\left(S,+\right)$ and $\left(S,\circ\right)$ both are inverse semigroups and they hold
\begin{align*}
    a\circ\left(b + c\right)
    = a\circ b - a + a\circ c \qquad \text{and} \qquad  a\circ a^-
    = - a + a,
\end{align*}
for all $a,b,c \in S$, where $-a$ and $a^-$ denote the inverses of $a$ with respect to $+$ and $\circ$, respectively, and we assume that the multiplication has higher precedence than the addition. As demonstrated in \cite[Thereom 10]{CaMaMiSt21x}, the additive semigroup of any weak brace is a Clifford semigroup. Moreover, if both $(S,+)$ and $(S,\circ)$ are groups, then $(S,+, \circ)$ is a \emph{skew brace}; in addition, if $(S,+)$ is abelian, then $(S,+, \circ )$ is a \emph{brace} \cite{Ru05}. Furthermore, weak braces form a subclass of \emph{inverse (left) semi-braces} \cite{CaMaSt21}, namely, triples $(S,+, \circ)$ with $(S,+)$ a semigroup and $(S, \circ)$ an inverse semigroup such that
$$a \circ (b+c)=a \circ b+a \circ \left(a^-+b\right),$$ for all $a,b \in S$. More specifically, weak braces are all the inverse semi-braces having $(S,+)$ as a Clifford semigroup and in which $-a+a\circ b=a \circ \left(a^-+b\right)$, for all $a,b \in S$. 
According to \cite[Theorem 11]{CaMaMiSt21x}, given a weak brace $(S,+, \circ)$, the map $r:S\times S\to S\times S$ defined by
 \begin{align*}
    r\left(a,b\right)
    = \left(-a+a\circ b, \left(-a+a\circ b\right)^-\circ a \circ b\right), 
 \end{align*}
for all $a,b\in S$, is a solution. Such a map $r$ has a behaviour close to bijectivity, since there exists the solution $r^{op}$ associated to the \emph{opposite weak brace}
$\left(S,+^{op}, \circ\right)$ of $S$, where $a+^{op}b:= b + a$, for all $a,b \in S$, such that
  \begin{align*}
      r\, r^{op}\, r = r, \qquad
      r^{op}\, r\, r^{op} = r^{op}, \qquad \text{and}\qquad rr^{op} = r^{op}r,
  \end{align*}
namely, $r$ is a
completely regular element in $\Map(S \times S)$. In this context, one can note that weak braces are instances of semitrusses that produces solutions but they are not necessarily  YB--semitrusses, see \cite{Br18} and \cite{CoJeVaVe21x}.

The aim of this paper is to show how to obtain weak braces from Rota--Baxter operators defined on Clifford semigroups. For this reason, we introduce the notion of  
 \emph{Rota--Baxter operator on a Clifford semigroup $(S,+)$}, i.e., 
a map $\mathfrak{R}:S \to S$ satisfying the relations
\begin{align*}\mathfrak{R}\left(a\right)+\mathfrak{R}\left(b\right)=\mathfrak{R}\left(a+\mathfrak{R}\left(a\right)+b-\mathfrak{R}\left(a\right)\right) \qquad \text{and} \qquad
    a+\mathfrak{R}\left(a\right)-\mathfrak{R}\left(a\right)=a,
\end{align*}
for all $a,b \in S$. Clearly, such definition includes the original one given for groups, cf. \cite[Definition 2.2]{GuLaSh20}.
Any Rota--Baxter operator $\mathfrak{R}$ on a Clifford semigroup $(S,+)$ determines a dual weak brace $S_\mathfrak{R}=(S, +, \circ)$ with
\begin{align*}
    a\circ b= a + \mathfrak{R}\left(a\right) + b - \mathfrak{R}\left(a\right),
\end{align*}
for all $a,b \in S$, that belongs to the class of weak braces having also $(S, \circ)$ as a Clifford semigroup. We name such particular weak braces \emph{dual weak braces}. The solution associated to a dual weak brace has a behaviour close to the non-degeneracy since also the maps $\lambda_a$ and $\rho_b$ are completely regular in
$\Map(S)$. In addition, we prove that the opposite solution $r^{op}$ of the solution $r$ associated to a dual weak brace $S_\mathfrak{R}$ derives from the dual weak brace obtained through the Rota--Baxter operator $\mathfrak{R}^{op}=-a+\mathfrak{R}(a)$, for every $a \in S$.\\
This fact has motivated us to deepen the structure of the dual weak brace on which we study the notion of ideal that allows us to define its quotient structure. In particular, we specialise ideals of dual weak braces obtained from Rota--Baxter operators. Besides, we introduce the \emph{socle} of a dual weak brace, a special ideal already investigated for braces \cite{Ru05} and other connected structures \cite{GuVe17, CaCoSt17, CeSmVe19,  CaCeSt20x}. Moreover, we show how to construct new ideals starting from two given ones.\\
Furthermore, we provide a method to construct Rota--Baxter operators on a Clifford semigroup seen as strong semilattices of groups, each of which admits a Rota--Baxter operator. Besides, we give a description of all endomorphisms $\mathfrak{R}$ of Clifford semigroups that are Rota--Baxter operators for which $\im \mathfrak{R}$ is commutative. Consistently with \cite[Proposition 5.1]{Ko21}, these special operators determine \emph{bi-weak braces}, i.e., special dual weak braces in which the roles of the sum and the multiplication can be reversed. Finally,  we characterise all Rota--Baxter endomorphisms of groups that are also idempotent.

\bigskip

\section{Basic results on weak braces}
\noindent The aim of this preliminary section is to give essential results on the structures of weak braces, recently introduced in \cite{CaMaMiSt21x} to find solutions of the Yang--Baxter equation. To this purpose, we initially recall some basics on Clifford semigroups, covered in detail in the books \cite{ClPr61, Ho95, Pe84}.
\medskip

A semigroup $S$ is called an \emph{inverse semigroup} if, for each $a\in S$, there exists a unique element $a^{-1}$ of $S$ such that $a=aa^{-1}a$ and $a^{-1}=a^{-1}aa^{-1}$. We call such an element $a^{-1}$ the \emph{inverse} of $a$. Evidently, every group is an inverse semigroup. The behavior  of inverse elements in an inverse semigroup $S$ is similar to that in a group, since $(a b)^{-1}=b^{-1} a^{-1}$ and $(a^{-1})^{-1}=a$, for all $a,b \in S$.\\
The set $\E(S)$ of the idempotents of an inverse semigroup $S$ is a commutative subsemigroup of $S$. Clearly, $e =e^{-1}$, for every $e \in \E(S)$.  In addition, the idempotents of $S$ are exactly the elements $a a^{-1}$ and $a^{-1}a$, with $a\in S$. 

If $S$ is an inverse semigroup such that $a a^{-1}=a^{-1}a$, for any $a\in S$, then $S$ is called \emph{Clifford semigroup}. In this case, we briefly write
\begin{align*}
   a^0:= a a^{-1}=a^{-1} a,
\end{align*}
for every $a\in S$. 
Equivalently, Clifford semigroups are the inverse semigroups in which the idempotents are central.\\
Moreover, if $\phi: S\to T$ is a homomorphism from a Clifford semigroup $S$ into a semigroup $T$, then $\im\phi$ is a Clifford semigroup, $\phi(a^{-1})=\phi(a)^{-1}$, and $\phi(a^{0})=\phi(a)^{0}$, for any $a\in S$ (cf. \cite[Lemma 2.4.4]{Ho95}).\\
According to \cite[Theorem 4.2.1]{Ho95}, an inverse semigroup $S$ is a Clifford semigroup if and only if it is a strong (lower) semilattice $Y$ of disjoint groups $G_{\alpha}$, for every $\alpha \in Y$.  Specifically, for each pair $\alpha,\beta$ of elements of $Y$ such that $\alpha \geq \beta$, let $\phii{\alpha}{\beta}:G_{\alpha}\to G_{\beta}$ be a homomorphism of groups such that
 \begin{enumerate}
     \item $\phii{\alpha}{\alpha}$ is the identical automorphism of $G_{\alpha}$, for every $\alpha \in Y$,
     \item $\phii{\beta}{\gamma}{}\phii{\alpha}{\beta}{} = \phii{\alpha}{\gamma}{}$, for all $\alpha, \beta, \gamma \in Y$ such that $\alpha \geq \beta \geq \gamma$.
 \end{enumerate} 
Then, $S = \displaystyle \bigcup_{\alpha\in Y}\,G_{\alpha}$ endowed with the operation defined by
\begin{align*}
    a\,b:= \phii{\alpha}{\alpha\beta}(a)\phii{\beta}{\alpha\beta}(b),
\end{align*}
for every $a\in G_{\alpha}$ and $b\in G_{\beta}$, is a \emph{strong semilattice $Y$ of groups $G_{\alpha}$}, usually written as $S=[Y; G_{\alpha}; \phii{\alpha}{\beta}]$.

\medskip

\begin{defin}\cite[Definition 5]{CaMaMiSt21x}\label{def:skew-inverse}
    Let $S$ be a set endowed with two operations $+$ and $\circ$ such that $\left(S,+\right)$ and $\left(S,\circ\right)$ are inverse semigroups. Then, $(S, +, \circ)$ is a \emph{weak (left) brace} if the following relations
 \begin{align*}
    a\circ\left(b + c\right)
    = a\circ b - a + a\circ c\qquad \text{and}\qquad a\circ a^-
    = - a + a
 \end{align*}
 hold, for all $a,b,c\in S$, where $-a$ and $a^-$ denote the inverses of $a$ with respect to $+$ and $\circ$, respectively. We call $(S,+)$ and $(S,\circ )$ the \emph{additive semigroup} and the \emph{multiplicative semigroup} of $\left(S,+,\circ\right)$, respectively, and, since the sets of idempotents $E(S,+)$ and $E(S, \circ)$ coincide, we denote it simply by $E(S)$.
\end{defin}

Given a weak brace $(S,+, \circ)$, let $\lambda_a:S\to S$ be the map defined by 
\begin{align*}
    \lambda_a\left(b\right) = - a + a\circ b,
\end{align*}
for all $a,b\in S$.  The map $\lambda_a$ can be also written as it is usual in the context of semi-braces \cite{CaCoSt17, JeAr19}, namely, $\lambda_a\left(b\right) 
    =  a\circ\left(a^- + b\right)$, for all $a,b\in S$.
Furthermore, the map $\lambda:S\to \End\left(S,+\right), a\mapsto \lambda_a$ is a homomorphism of the inverse semigroup $\left(S,\circ\right)$ into the endomorphism semigroup of $\left(S,+\right)$.\\
The operations of sum and multiplication of any weak brace are closely related by the following relation
\begin{align}\label{eq:circ+}
    a\circ b = a + \lambda_a\left(b\right),
\end{align}
for all $a,b\in S$, as proved in \cite[Lemma 1]{CaMaMiSt21x}. In particular, they coincide on $E(S)$. Indeed, if $e_1,e_2\in E(S)$, we get
\begin{align*}
    e_1+ e_2 &= \left(e_1+e_2\right)\circ \left(e_1+e_2\right)
    = \left(e_1+e_2\right) \circ e_1-\left(e_1+e_2\right)+\left(e_1+e_2\right) \circ e_2\\
    &= e_1 \circ \left(e_1+e_2\right)+e_1+e_2+e_2 \circ \left(e_2+e_1\right) =e_1+\lambda_{e_1}\left(e_2\right)+e_2+\lambda_{e_2}\left(e_1\right)\\
    &=e_1\circ e_2+e_2 \circ e_1 &\mbox{by \eqref{eq:circ+}} \\
    &=e_1\circ e_2.
\end{align*}

As proved in \cite[Theorem 8]{CaMaMiSt21x}, the additive semigroup of any weak brace $S$ is a Clifford semigroup. 
Moreover, the multiplicative semigroup of a weak brace need not be a Clifford semigroup (see, for instance, \cite[Example 3]{CaMaMiSt21x}).
Thus, the following question arises.
\begin{que} Describe the structure of the multiplicative semigroup of a weak brace.
\end{que}

\noindent However, one can find some instances of weak braces having also a Clifford semigroup as a multiplicative structure among  \cite{CaCoSt21-gen, CaMaMiSt21x, CaMaSt21}.
\begin{defin}
A  weak brace $(S, +,\circ)$ is said to be a \emph{dual weak brace} if $\left(S,\circ\right)$ is a Clifford semigroup.
\end{defin}
\noindent In the case of a dual weak brace $S$, we highlight that $a^0$ is also $a^-\circ a$, for every $a \in S$.\\
Note that any Clifford semigroup $\left(S, \circ\right)$ determines two dual weak braces, by setting $a + b:= a\circ b$ or $a + b:= b\circ a$ for all $a,b\in S$. The first weak brace is the \emph{trivial weak brace} and the second one the \emph{almost trivial weak brace} (see \cite[Example 1]{CaMaMiSt21x}).
\medskip

Below, we show that in any weak brace the sum can be written in terms of the multiplication as it is already known in the context of skew braces, see \cite[Remark 1.8]{GuVe17}.
\begin{prop}\label{+circ}
Let $S$ be a weak brace. Then, 
$a+b=a \circ \lambda_{a^-}\left(b\right),$
for all $a,b\in S$.
\begin{proof}
Let $a,b \in S$. Observing that by \cite[Lemma 1]{CaMaMiSt21x} $a \circ (-b)=a-a\circ b+a$, we obtain
\begin{align*}
\left(a\circ\lambda_{a^-}\left(b\right)\right)^-
&= \left(a+b\right)^-\circ\left(a - a\right)
= \left(a+b\right)^-\circ a - \left(a+b\right)^- + \left(a+b\right)^-\circ\left(-a\right)\\
&= \left(a+b\right)^-\circ a  - \left(a+b\right)^-\circ a + \left(a+b\right)^-\\
&= \left(a+b\right)^- - \left(a+b\right)^-\circ a  +\left(a+b\right)^-\circ a\\
&=  \left(a+b\right)^-  + \lambda_{ \left(a+b\right)^-}\left(-a\right)  -  \left(a+b\right)^- +  \left(a+b\right)^-\circ a &\mbox{by \eqref{eq:circ+}}\\
&= \left(a+b\right)^-  +  \left(a+b\right)^-\circ  \left(a + b  - a + a\right)\\
&= \left(a+b\right)^-  +  \left(a+b\right)^-\circ  \left(a + b\right)\\
&= \left(a+b\right)^-  -  \left(a+b\right)^- + \left(a + b\right)^-\\
&= \left(a+b\right)^-.
\end{align*}
Therefore, the claim follows.
\end{proof}
\end{prop}

\medskip

The trivial and the almost trivial weak braces lie in a special class of dual weak braces that we call \emph{bi-weak braces}.
We underline that such a definition has been originally introduced by Childs in \cite{Ch19} for the particular case of skew braces. 
\begin{defin}
Let $(S, +, \circ)$ be a weak brace. Then, $S$ is called a \emph{bi-weak brace} if $(S, \circ, +)$ is a weak brace, i.e., the relations
\begin{align*}
    a + ( b \circ c)= \left(a +b\right) \circ a^- \circ \left(a +c\right)
    \qquad\text{and}\qquad
    a - a = a^-\circ a
\end{align*}
hold, for all $a,b,c \in S$.
\end{defin}

\medskip

The following is a general way to construct weak braces  involving strong semilattices of inverse semi-braces provided in \cite[Proposition 4]{CaMaSt21} and also used in the case of generalized left semi-braces in \cite{CaCoSt21-gen}. 

\begin{prop}\label{prop:Strong-Lattice-Inverse-Semi-Brace}
Let $Y$ be a (lower) semilattice, $\left\{S_{\alpha}\ \left|\ \alpha \in Y\right.\right\}$ a family of disjoint weak braces. For each pair $\alpha,\beta$ of elements of $Y$ such that $\alpha \geq \beta$, let $\phii{\alpha}{\beta}:\s{\alpha}\to \s{\beta}$ be a homomorphism of weak braces such that
\begin{enumerate}
    \item $\phii{\alpha}{\alpha}$ is the identical automorphism of $\s{\alpha}$, for every $\alpha \in Y$,
    \item $\phii{\beta}{\gamma}{}\phii{\alpha}{\beta}{} = \phii{\alpha}{\gamma}{}$, for all $\alpha, \beta, \gamma \in Y$ such that $\alpha \geq \beta \geq \gamma$.
\end{enumerate}
Then, $S = \bigcup\left\{\s{\alpha}\ \left|\ \alpha\in Y\right.\right\}$ endowed with the addition and the multiplication defined by
\begin{align*}
    a+b:= \phii{\alpha}{\alpha\beta}(a)+\phii{\beta}{\alpha\beta}(b)
    \quad\text{and}\quad
     a\circ b:= \phii{\alpha}{\alpha\beta}(a)\circ\phii{\beta}{\alpha\beta}(b),
\end{align*}
for all $a\in \s{\alpha}$ and $b\in \s{\beta}$, is a weak brace, called \emph{strong semilattice $S$ of weak braces $S_{\alpha}$} and is denoted by $S=[Y; S_\alpha;\phii{\alpha}{\beta}]$.
\begin{proof}
Since $S$ is an inverse semi-brace and $(S,+)$ clearly is an inverse semigroup, by \cite[Proposition 16]{CaMaMiSt21x},  it is sufficient to show that  $a\circ\left(a^-+b\right)=-a+a\circ b$, for all $a,b \in S$.
If $a\in \s{\alpha}$ and $b\in \s{\beta}$, then
\begin{align*}
    a\circ\left(a^- + b\right)
    &=  a\circ\left(\phii{\alpha}{\alpha\beta}\left(a^-\right) + \phii{\beta}{\alpha\beta}\left(b\right)\right)\\
    &= \phii{\alpha}{\alpha\beta}\left(a\right)\circ
    \phii{\alpha\beta}{\alpha\beta}\left(\phii{\alpha}{\alpha\beta}\left(a^-\right) + \phii{\beta}{\alpha\beta}\left(b\right)\right)\\
    &=\phii{\alpha}{\alpha\beta}\left(a\right)\circ \left(\phii{\alpha}{\alpha\beta}\left(a\right)^- + \phii{\beta}{\alpha\beta}\left(b\right)\right)\\ 
    &= -\phii{\alpha}{\alpha\beta}\left(a\right)
    +\phii{\alpha}{\alpha\beta}\left(a\right)\circ 
    \phii{\beta}{\alpha\beta}\left(b\right) &\mbox{$S_{\alpha\beta}$ is a weak brace} \\
    &= \phii{\alpha}{\alpha\beta}\left(-a\right)
    +\phii{\alpha\beta}{\alpha\beta}\left(\phii{\alpha}{\alpha\beta}\left(a\right)\circ 
    \phii{\beta}{\alpha\beta}\left(b\right)\right)\\
    &= - a + \phii{\alpha}{\alpha\beta}\left(a\right)\circ \phii{\beta}{\alpha\beta}\left(b\right)\\
    &= - a + a\circ b.
\end{align*}
Therefore, the claim follows.
\end{proof}
\end{prop}
\medskip 

\begin{cor}
A strong semilattice $S=[Y; S_\alpha;\phii{\alpha}{\beta}]$ is a dual weak brace if and only if $S_{\alpha}$ is a dual weak brace, for every $\alpha\in Y$, as well. 
Specifically, $S$ is a bi-weak brace if and only if $S_{\alpha}$ is a bi-weak brace.
\end{cor}

\bigskip

\section{Rota--Baxter operators on Clifford semigroups}

\noindent In this section, we present and investigate the Rota--Baxter operators on a Clifford semigroup, consistently with that introduced for groups in \cite{GuLaSh20}. Moreover, we provide some examples.

\medskip

\begin{defin}\label{def:Rb-oper}
Let $(S,+)$ be a Clifford semigroup. A map $\mathfrak{R}:S \to S$ satisfying
\begin{align*}\mathfrak{R}\left(a\right)+\mathfrak{R}\left(b\right)=\mathfrak{R}\left(a+\mathfrak{R}\left(a\right)+b-\mathfrak{R}\left(a\right)\right) \qquad \text{and} \qquad
    a+\mathfrak{R}\left(a\right)^0=a,
\end{align*}
for all $a,b \in S$, is called \emph{Rota--Baxter operator} on $(S,+)$.
\end{defin}
\noindent  If $(S,+)$ and $(T,+)$ are Clifford semigroups, then two Rota--Baxter operators $\mathfrak{R}$ and $\mathfrak{T}$ on $S$ and $T$, respectively, are said to be \emph{equivalent} if there exists an isomorphism $\varphi$ from $(S,+)$ to $(T,+)$ such that $\mathfrak{R}=\varphi^{-1}\mathfrak{T}\varphi$.

\medskip

\begin{ex}
Let $\left(S, +\right)$ be a Clifford semigroup and $\varphi$ an idempotent endomorphism of $\left(S, +\right)$ such that $\varphi\left(e\right) = e$, for every $e\in E\left(S,+\right)$. Then, the map $\mathfrak{R}:= -\varphi$ is a Rota--Baxter operator on $S$. 
\end{ex}
\noindent Easy examples of operators belonging to the previous class are the \emph{elementary Rota--Baxter operators}, namely, the
maps $\mathfrak{E}$ and $\mathfrak{O}$ from a Clifford semigroup $(S,+)$ into itself given by $\mathfrak{E}(a) = a^0$ and $\mathfrak{O}(a) = -a$, for every $a \in S$. 
Such a terminology is analogue to that used in \cite{BaGu21groupsx}. One can observe that in the context of groups such Rota--Baxter operators are not equivalent.
In addition, if $(S,+)$ is a Clifford monoid with identity $0$, the map $\mathfrak{R}(a) = 0$, for every $a\in S$, is a Rota--Baxter operator on $(S,+)$. 

If $\mathfrak{R}$ is a Rota--Baxter operator on a Clifford semigroup $(S,+)$, the map $\mathfrak{R}^{op}:S \to S$, given by $\mathfrak{R}^{op}(a)=-a+\mathfrak{R}(-a)$, for every $a \in S$, is still a Rota--Baxter operator on $(S,+)$, which we call \emph{opposite Rota--Baxter operator} of $\mathfrak{R}$. We observe that the elementary Rota--Baxter operators can be obtained one from the other since $\mathfrak{O}=\mathfrak{E}^{op}$ and $\mathfrak{E}=\mathfrak{O}^{op}$.

\medskip

The following are two examples of Rota--Baxter operators on a Clifford monoid which admits an exact factorization. We recall that a Clifford semigroup $(S, +)$ admits an \emph{exact factorization} $(U,V)$ if $U$ and $V$ are Clifford subsemigroups of $S$ and any element $a\in S$ can be written in a unique way as $a = u_a + v_a$, with $u_a\in U$ and $v_a\in V$. In such a case, it holds that $U\cap V= \{0\}$, with $0$ both right identity of $U$ and left identity of $V$, cf. \cite[Theorem 2]{Ca87}. As a consequence, $(S, +)$ is a Clifford monoid.

\begin{ex}
   Let $(U, V)$ be an exact factorization of a Clifford monoid $(S,+)$.
 Then, the maps $\mathfrak{R}:S\to S$ and $\mathfrak{T}:S\to S$ given by
 \begin{align*}
     \mathfrak{R}(a)= - v_a \qquad \text{and} \qquad \mathfrak{T}(a)= u_a^0- v_a
 \end{align*}
 for every $a \in S$, respectively, are two Rota--Baxter operators on $S$. \\
We underline that the map $\mathfrak{R}$ is a Rota--Baxter operator also in the case of a right univocal factorizable Clifford semigroup. To this purpose, we recall that, in general, a factorization $\left(U,V\right)$ of a semigroup $S$ is \emph{right univocal} if 
\begin{align*}
    \forall \ u,u'\in U, \ v,v'\in V\quad u + v = u' + v'\implies v = v',
\end{align*}
see, for instance, \cite[Definition 1]{Ca87}.
\end{ex}

\bigskip

In the next result, we introduce a method for constructing a new Rota--Baxter operator on a Clifford semigroup seen as strong semilattice $Y$ of groups $G_{\alpha}$ starting from given Rota--Baxter operators on each group $G_{\alpha}$. Hereinafter, we denote by $0_{\alpha}$ the identity of $G_{\alpha}$, for every $\alpha \in Y$.

\begin{theor}\label{teo_RB_semilattice}
Let $(S,+)$ be a Clifford semigroup, $S=[Y; G_{\alpha}; \phii{\alpha}{\beta}]$ and assume that $\mathfrak{R}_{\alpha}$ is a Rota--Baxter operator on the group $(G_{\alpha},+)$, for every $\alpha \in Y$. 
Then, the map $\mathfrak{R}:S \to S$ given by
$\mathfrak{R}\left(a\right) = \mathfrak{R}_{\alpha}\left(a\right)$, for every $a\in G_{\alpha}$, is a Rota--Baxter operator on $(S,+)$ if and only if the condition
\begin{align}\label{condizione_semireticolo}
    \mathfrak{R}_{\beta}\phii{\alpha}{\beta}
    = \phii{\alpha}{\beta}\mathfrak{R}_{\alpha},
\end{align}
is satisfied, for all $\alpha, \beta\in Y$ such that $\alpha\geq \beta$. We call such an operator $\mathfrak{R}$ \emph{strong Rota--Baxter operator on $(S,+)$}.
\begin{proof}
Initially, assume that \eqref{condizione_semireticolo} holds.
If $a\in G_{\alpha}$ and $b\in G_{\beta}$, we have that
\begin{align*}
    a + \mathfrak{R}\left(a\right) + b - \mathfrak{R}\left(a\right)
    &= \phii{\alpha}{\alpha\beta}\left(a + \mathfrak{R}_{\alpha}\left(a\right)\right) + \phii{\beta}{\alpha\beta}\left(b\right) - \mathfrak{R}_{\alpha}\left(a\right)\\
    &= \phii{\alpha}{\alpha\beta}\left(a\right) + \phii{\alpha}{\alpha\beta}\mathfrak{R}_{\alpha}\left(a\right)
    + \phii{\beta}{\alpha\beta}\left(b\right) 
    - \phii{\alpha}{\alpha\beta}\mathfrak{R}_{\alpha}\left(a\right)\\
    &= \phii{\alpha}{\alpha\beta}\left(a\right)+\mathfrak{R}_{\alpha\beta}\phii{\alpha}{ \alpha\beta}\left(a\right)+\phii{\beta}{ \alpha\beta}\left(b\right)-\mathfrak{R}_{\alpha\beta}\phii{\alpha}{ \alpha\beta}\left(a\right)&\mbox{by \eqref{condizione_semireticolo}}
\end{align*}
Thus, we obtain
\begin{align*}
    &\mathfrak{R}\left(a+\mathfrak{R}\left(a\right)+b-\mathfrak{R}\left(a\right)\right)\\
    &=\mathfrak{R}_{\alpha\beta}\left(\phii{\alpha}{ \alpha\beta}\left(a\right)+\mathfrak{R}_{\alpha\beta}\phii{\alpha}{ \alpha\beta}\left(a\right)+\phii{\beta}{\alpha\beta}\left(b\right)-\mathfrak{R}_{\alpha\beta}\phii{\alpha}{ \alpha\beta}\left(a\right)\right) \\
    &= \mathfrak{R}_{\alpha\beta}\phii{\alpha}{ \alpha\beta}\left(a\right)+\mathfrak{R}_{\alpha\beta}\phii{\beta}{\alpha\beta}\left(b\right)\\
    &=\phii{\alpha}{\alpha\beta}\mathfrak{R}_{\alpha} \left(a\right)+ \phii{\beta}{\alpha\beta}\mathfrak{R}_{\beta} \left(b\right)&\mbox{by \eqref{condizione_semireticolo}}\\
    &=\mathfrak{R}\left(a\right)+\mathfrak{R}\left(b\right).
\end{align*}
Therefore, since the second condition in \cref{def:Rb-oper} is trivially satisfied, $\mathfrak{R}$ is a Rota--Baxter operator on $S$.\\
Conversely, assume that $\mathfrak{R}$ is a Rota--Baxter operator on $(S,+)$. Then, if $a\in G_{\alpha}$, by \cite[Lemma 2.5]{BaGu21groupsx} we get 
\begin{align*}
    \mathfrak{R}\left(a\right) + \mathfrak{R}\left(0_{\beta}\right)
    = \mathfrak{R}_{\alpha}\left(a\right) + \mathfrak{R}_{\beta}\left(0_{\beta}\right)
    = \phii{\alpha}{\alpha\beta}\mathfrak{R}_{\alpha}\left(a\right) + 0_{\alpha\beta}
    =  \phii{\alpha}{\alpha\beta}\mathfrak{R}_{\alpha}\left(a\right).
\end{align*}
On the other hand,
\begin{align*}
    \mathfrak{R}\left(a + \mathfrak{R}\left(a\right) + 0_{\beta} - \mathfrak{R}\left(a\right)\right) 
    &= \mathfrak{R}_{\alpha\beta}\left(\phii{\alpha}{\alpha\beta}\left(a\right) + \phii{\alpha}{\alpha\beta}\mathfrak{R}_{\alpha}\left(a\right)
    + \phii{\beta}{\alpha\beta}\left(0_{\beta}\right) 
    - \phii{\beta}{\alpha\beta}\mathfrak{R}_{\alpha}\left(a\right)\right)\\
    &= \mathfrak{R}_{\alpha\beta}\phii{\alpha}{\alpha\beta}\left(a\right),
\end{align*}
therefore \eqref{condizione_semireticolo} holds.
\end{proof}
\end{theor}

We observe that not all Rota-Baxter operators on Clifford semigroups
can be constructed as in \cref{teo_RB_semilattice}. Indeed, if $(S,+)$ is a Clifford monoid with $|E(S)| > 1$, the  map $\mathfrak{R}:S\to S$ given by $\mathfrak{R}(a) = 0$, for every $a \in S$, is such an example. 
\medskip

\noindent Strong Rota--Baxter operators on a Clifford semigroup $S=[Y; G_{\alpha}; \phii{\alpha}{\beta}]$ can be easily obtained starting from arbitrary Rota--Baxter operators on groups $G_{\alpha}$.
\begin{ex}
Let $S=[ Y; G_{\alpha}; \phii{\alpha}{\beta}]$ be a Clifford semigroup with $ \phii{\alpha}{\beta}(a) = 0_{\beta}$, for every $a\in G_{\alpha}$ and $\alpha,\beta\in Y$ such that $\alpha\geq \beta$. If $\mathfrak{R}_{\alpha}$ is an arbitrary Rota--Baxter operator on $G_{\alpha}$, for every $\alpha \in Y$, then the map $\mathfrak{R}:S \to S$ given by $\mathfrak{R}(a) = \mathfrak{R}_{\alpha}(a)$, for every $a\in G_{\alpha}$, is a strong Rota--Baxter operator on $S$.
\end{ex}

\bigskip

Finally, let us show some useful properties of Rota--Baxter operators on Clifford semigroups.

\begin{prop}\label{prop:RB-prop-E}
Let $\mathfrak{R}$ be a Rota--Baxter operator on a Clifford semigroup $(S,+)$. Then, 
\begin{align*}
\mathfrak{R}\left(E\left(S\right)\right)\subseteq E\left(S\right).
\end{align*}
In particular, if $(S,+)$ is a Clifford monoid with identity $0$, it holds that $\mathfrak{R}\left(0\right)=0$.
\begin{proof}
If $e\in E\left(S\right)$, we have that
\begin{align*}
    \mathfrak{R}\left(e\right)
    =  \mathfrak{R}\left(e + \mathfrak{R}\left(e\right) -  \mathfrak{R}\left(e\right)\right)
    =  \mathfrak{R}\left(e +  \mathfrak{R}\left(e\right) + e -  \mathfrak{R}\left(e\right)\right)
    =  \mathfrak{R}\left(e\right) +  \mathfrak{R}\left(e\right).
\end{align*}
Besides , if $(S,+)$ is a Clifford monoid with identity $0$, since $\mathfrak{R}(0)$ is idempotent, we get
\begin{align*}
    0 = 0 + \mathfrak{R}\left(0\right)-\mathfrak{R}\left(0\right)
    = \mathfrak{R}\left(0\right)+\mathfrak{R}\left(0\right)
    = \mathfrak{R}\left(0\right),
\end{align*}
thus the claim follows.
\end{proof}
\end{prop}

In general, it is not true that $E(S) \subseteq \mathfrak{R}\left(E\left(S\right)\right)$. In fact, it is enough to take the Rota-Baxter operator of constant value $0$ on a Clifford monoid $(S,+)$ with $|E(S)|>1$.

\medskip

\begin{prop}\label{prop:RB-prop}
Let $\mathfrak{R}$ be a Rota--Baxter operator on a Clifford semigroup $(S,+)$.
Then, the following statements hold:
\begin{enumerate}
    \item $\mathfrak{R}(a)=\mathfrak{R}(a) + \mathfrak{R}\left(a^0\right)$,
    \item $a= a+\mathfrak{R}\left(a^0\right)
    $,
    \item $- \mathfrak{R}\left(a\right) = \mathfrak{R}\left(-\mathfrak{R}\left(a\right) - a + \mathfrak{R}\left(a\right)\right) $,
\end{enumerate}
for every $a\in S$.
\begin{proof}
If $a\in S$, we have that
\begin{align*}
\mathfrak{R}(a) 
= \mathfrak{R}\left(a + \mathfrak{R}\left(a\right) ^0\right)
= \mathfrak{R}\left(a + \mathfrak{R}(a) + a^0 - \mathfrak{R}(a)\right) 
= \mathfrak{R}\left(a\right) + \mathfrak{R}\left(a^0\right),
\end{align*}
i.e., $1.$ is verified. Moreover, by $1.$,
\begin{align*}
   a=a - \mathfrak{R}(a)+ \mathfrak{R}(a)
  =a - \mathfrak{R}(a)+ \mathfrak{R}(a)+\mathfrak{R}\left(a^0\right) 
   = 
   a + \mathfrak{R}\left(a^0\right),
\end{align*}
and so $2.$ holds. To prove $3.$, we note that
\begin{align}\label{relazione}
    \mathfrak{R}(a) +\mathfrak{R}\left(-\mathfrak{R}\left(a\right) - a + \mathfrak{R}\left(a\right)\right)=\mathfrak{R}\left(a^0\right).
\end{align}
Thus, we obtain that
\begin{align*}
    \mathfrak{R}(a) +\mathfrak{R}\left(-\mathfrak{R}\left(a\right) - a + \mathfrak{R}\left(a\right)\right)+\mathfrak{R}(a)&=\mathfrak{R}\left(a^0\right)+\mathfrak{R}(a) &\mbox{by \eqref{relazione}}\\
    &=\mathfrak{R}(a) + \mathfrak{R}\left(a^0\right)&\mbox{by \cref{prop:RB-prop-E}}\\
    &= \mathfrak{R}(a)&\mbox{by $1.$}
\end{align*} 
and
\begin{align*}
&\mathfrak{R}\left(-\mathfrak{R}\left(a\right) - a + \mathfrak{R}\left(a\right)\right)+ \mathfrak{R}(a)+  \mathfrak{R}\left(-\mathfrak{R}\left(a\right) - a + \mathfrak{R}\left(a\right)\right)  \\
&\quad=\mathfrak{R}\left(-\mathfrak{R}\left(a\right) - a + \mathfrak{R}\left(a\right)\right)+  \mathfrak{R}\left(a^0\right) &\mbox{by \eqref{relazione}}\\
&\quad=\mathfrak{R}\left(a^0\right)+ \mathfrak{R}\left(-\mathfrak{R}\left(a\right) - a + \mathfrak{R}\left(a\right)\right) &\mbox{by \cref{prop:RB-prop-E}}\\
&\quad=\mathfrak{R}\left(a^0+\mathfrak{R}\left(a^0\right)-\mathfrak{R}\left(a\right) - a + \mathfrak{R}\left(a\right)-\mathfrak{R}\left(a^0\right)\right)\\
&\quad=\mathfrak{R}\left(-\mathfrak{R}\left(a\right) - a + \mathfrak{R}\left(a\right)+\mathfrak{R}\left(a^0\right)\right)&\mbox{by \cref{prop:RB-prop-E}}\\
&\quad=\mathfrak{R}\left(-\mathfrak{R}\left(a\right) - a + \mathfrak{R}\left(a\right)\right)  &\mbox{by $1.$}
\end{align*}
hence $3.$ follows. 
\end{proof}
\end{prop}

\medskip

\begin{rem}\label{remmeno}
If $\mathfrak{R}$ is a Rota--Baxter operator on a Clifford semigroup $(S,+)$, by $3.$ in \cref{prop:RB-prop}, we obtain that $\im\mathfrak{R}$ is a Clifford subsemigroup of $(S,+)$.
Furthermore, we can make explicit the statement in \cref{prop:RB-prop-E}, since, using $3.$ in \cref{prop:RB-prop} and \eqref{relazione}, we get
\begin{align}\label{eq:rem:ImR}
 \mathfrak{R}\left(a^0\right)
    = \mathfrak{R}\left(a\right)^0,  
\end{align}
for every $a\in S$.
\end{rem}

\bigskip

\section{Rota--Baxter endomorphisms}

 \noindent This section focuses on endomorphisms of Clifford semigroups that are Rota--Baxter operators, namely Rota--Baxter endomorphisms. Specifically, we give a description of such maps $f$ for which $\im f$ is commutative. In addition,  we characterize Rota--Baxter endomorphisms of groups that are also idempotent.

\medskip

\begin{defin}
Let $(S,+)$ be a Clifford semigroup. Any endomorphism  $\mathfrak{R}: S \to S$ that is a Rota--Baxter operator is called a \emph{Rota--Baxter endomorphism} of $(S,+)$.
\end{defin}
\noindent Clearly, the constant map of value $0$ of a Clifford monoid $(S,+)$ is an easy example of Rota--Baxter endomorphism.

\medskip

From now on, given two Clifford semigroups $(S,+)$ and $(T,+)$, we will denote by $\Hom\left(S,T\right)$ the set of all the homomorphisms from $(S,+)$ to $(T,+)$. If  $f\in\Hom\left(S,T\right)$ is such that $\im f$ is commutative, we say that $f$ is a \emph{commutative homomorphism}. In the context of groups,  these maps are named \emph{abelian homomorphisms} \cite{Ch12}.\\ 
If $G$ is a group, one has that the abelian endomorphisms of $G$ are Rota--Baxter operators.
Instead, if $\left(S, +\right)$ is a Clifford semigroup, not every commutative endomorphism of $(S,+)$ is a Rota--Baxter operator. Indeed, if $|E(S)| > 1$ and $e\in E(S)$, the constant map of value $e$ is such an example. Thus, our aim is to show how to construct all the commutative Rota--Baxter endomorphisms on a Clifford semigroup $S = [Y; G_{\alpha}; \phii{\alpha}{\beta}]$ starting from given abelian homomorphisms. To this purpose, we recall the following result.
\begin{lemma}[cf. Proposition II.2.8, \cite{Pe84}]\label{lemma:Worawiset}
Let $S = [Y; G_{\alpha}; \phii{\alpha}{\beta}]$ be a Clifford semigroup, $h:Y \to Y$ an endomorphism, and consider a family of maps $f_{\alpha}\in \Hom\left(G_{\alpha}, G_{h\left(\alpha\right)}\right)$, for every $\alpha\in Y$, such that
\begin{align*}
    f_{\beta}\phi_{\alpha, \beta}=\phi_{h\left(\alpha\right), h\left(\beta\right)}f_{\alpha}
\end{align*}
holds, for all $\alpha, \beta \in Y$ such that $\alpha \geq\beta$. Then, the map $f: S \to S$ given by $f(a)=f_{\alpha}(a)$, for every $a \in G_{\alpha}$, is an endomorphism of $S$. 
Conversely, every endomorphism of $S$ can be so constructed.
\end{lemma}

\medskip
\begin{theor}\label{th:comm_RB-end}
Let $S = [Y; G_{\alpha}; \phii{\alpha}{\beta}]$ be a Clifford semigroup, $h:Y \to Y$ an endomorphism such that $\alpha\leq h\left(\alpha\right)$, and consider a family of abelian homomorphisms $f_{\alpha}\in \Hom\left(G_{\alpha}, G_{h\left(\alpha\right)}\right)$, for every $\alpha\in Y$. 
If the following condition is satisfied
\begin{align}\label{eq:Worawiset}
    f_{\beta}\phi_{\alpha, \beta}=\phi_{h\left(\alpha\right), h\left(\beta\right)}f_{\alpha},
\end{align}
for all $\alpha, \beta \in Y$ such that $\alpha \geq\beta$, then the map $\mathfrak{R}: S \to S$ given by $\mathfrak{R}(a)=f_{\alpha}(a)$, for every $a \in G_{\alpha}$, is a commutative Rota--Baxter endomorphism of $S$.
Conversely, any commutative Rota--Baxter endomorphism on $S$ can be so constructed.
\begin{proof}
Initially, by \cref{lemma:Worawiset}, note that $\mathfrak{R}$ is an endomorphism of $S$ and obviously $\mathfrak{R}$ is commutative since $f_{\alpha}$ is abelian, for every $\alpha\in Y$.
Now, if  $a\in G_{\alpha}$ and $b\in G_{\beta}$, with $\alpha,\beta\in Y$ such that $\alpha \geq \beta$, we have that
\begin{align*}
a&+ \mathfrak{R}\left(a\right)+b-\mathfrak{R}\left(a\right)\\
&= \left(\phi_{\alpha,\alpha h\left(\alpha\right)}\left(a\right) 
+ \phi_{h\left(\alpha\right),\alpha h\left(\alpha\right)}f_{\alpha}\left(a\right)\right)
+ \left(\phi_{\beta,\beta h\left(\alpha\right)}\left(b\right)
- \phi_{h\left(\alpha\right),\beta h\left(\alpha\right)}f_{\alpha}\left(a\right)
\right)\\
&= a + \phi_{h\left(\alpha\right),\alpha}f_{\alpha}\left(a\right) 
+ \phi_{\beta,\beta h\left(\alpha\right)}\left(b\right)
- \phi_{h\left(\alpha\right), \beta h\left(\alpha\right)}f_{\alpha}\left(a\right)\\
&=\phi_{\alpha,\alpha\beta h\left(\alpha\right)}
\left(a + \phi_{h\left(\alpha\right),\alpha}f_{\alpha}\left(a\right)\right) 
+ \phi_{\beta h\left(\alpha\right),\alpha\beta}\left(\phi_{\beta,\beta h\left(\alpha\right)}\left(b\right)
- \phi_{h\left(\alpha\right), \beta h\left(\alpha\right)}f_{\alpha}\left(a\right)\right)\\
&= \phi_{\alpha,\alpha\beta}\left(a\right)
+ \phi_{\alpha,\alpha\beta}\phi_{h\left(\alpha\right),\alpha}f_{\alpha}\left(a\right)
+ \phi_{\beta,\alpha\beta}\left(b\right)
- \phi_{h\left(\alpha\right),\alpha\beta}f_{\alpha}\left(a\right)\\
&=\phi_{\alpha, \alpha\beta}\left(a\right)+\phi_{ h\left(\alpha\right), \alpha\beta}f_{\alpha}\left(a\right)+\phi_{\beta, \alpha\beta}\left(b\right)-\phi_{ h\left(\alpha\right), \alpha\beta}f_{\alpha}\left(a\right)
\end{align*}
and, consequently,
\begin{align*}
\mathfrak{R}&\left(a+ \mathfrak{R}\left(a\right)+b-\mathfrak{R}\left(a\right)\right)\\
&= f_{\alpha\beta}\left(\phi_{\alpha, \alpha\beta}\left(a\right)
   + \phi_{h\left(\alpha\right), \alpha\beta}f_{\alpha}\left(a\right) 
   + \phi_{\beta,\alpha\beta}\left(a\right)
   - \phi_{h\left(\alpha\right), \alpha\beta}f_{\alpha}\left(a\right)\right)\\
   &=f_{\alpha\beta}\phii{\alpha}{\alpha\beta}\left(a\right) + f_{\alpha\beta}\phii{\beta}{\alpha\beta}\left(b\right)
   &\mbox{$f_{\alpha\beta}$ is abelian}\\
   &= \phii{h\left(\alpha\right)}{h\left(\alpha\beta\right)}f_{\alpha}\left(a\right)
   + \phii_{h\left(\beta\right), h\left(\alpha\beta\right)}f_{\beta}\left(b \right)
   &\text{by \eqref{eq:Worawiset}}\\
      &=  f_{\alpha}\left(a\right) + f_{\beta}\left(b\right)=\mathfrak{R}\left(a\right) + \mathfrak{R}\left(b\right)
\end{align*}
Moreover, it holds that
\begin{align*}
    a + \mathfrak{R}\left(a\right)^0
    &= a + f_{\alpha}\left(0_{\alpha}\right)
    = \phii{\alpha}{\alpha h\left(\alpha\right)}\left(a\right) + 0_{\alpha h\left(\alpha\right)}
    = \phii{\alpha}{\alpha}\left(a\right) + 0_{\alpha}
    = a.
\end{align*}
Therefore, $\mathfrak{R}$ is a commutative Rota--Baxter endomorphism of $S$.\\
Conversely, assume that $\mathfrak{R}$ is a commutative Rota--Baxter endomorphism of $S$. Let $h:Y\to Y$ be the map such that $h\left(\alpha\right)=\chi$ if $\mathfrak{R}\left(0_{\alpha}\right) = 0_{\chi}$. Note that $h\left(\alpha\right) \geq \alpha$, since if $a \in G_{\alpha}$, by $2.$ in \cref{prop:RB-prop}, we obtain
\begin{align*}
a=a+\mathfrak{R}\left(a^0\right)=a+0_{\chi} \in G_{\alpha\chi},
\end{align*}
thus $\alpha\chi=\alpha$, and so $h(\alpha)\geq \alpha$.
Now, denoted by $f_\alpha:= \mathfrak{R}_{|_{G_{\alpha}}}$, the map $\mathfrak{R}_{\alpha}:= \phii{h\left(\alpha\right)}{\alpha}f_{\alpha}$ is clearly an abelian endomorphism of $G_{\alpha}$, for every $\alpha\in Y$. Moreover, if $a\in G_{\alpha}$ and $\alpha,\beta\in Y$ with $\alpha\geq \beta$, we get
\begin{align*}
    \mathfrak{R}_{\beta}\phii{\alpha}{\beta}\left(a\right)
    &= \phii{h\left(\beta\right)}{\beta}f_{\beta}\phii{\alpha}{\beta}\left(a\right)
    = \phii{h\left(\beta\right)}{\beta}\mathfrak{R}\left(a + 0_{\beta}\right)
    = \phii{h\left(\beta\right)}{\beta}\left(\mathfrak{R}\left(a\right) + \mathfrak{R}\left(0_{\beta}\right)\right)\\
    &=\phii{h\left(\beta\right)}{\beta}\phii{h\left(\alpha\right)}{h\left(\alpha\beta\right)}\mathfrak{R}\left(a\right)
    = \phii{h\left(\alpha\right)}{\beta}f_{\alpha}\left(a\right)
    = \phii{\alpha}{\beta}\phii{h\left(\alpha\right)}{\alpha}f_{\alpha}\left(a\right)
    = \phii{\alpha}{\beta}R_{\alpha}\left(a\right).
\end{align*}
Therefore, the claim follows.
\end{proof}
\end{theor}

\medskip

\noindent Obviously, any Rota--Baxter endomorphism is a strong Rota--Baxter operator only in the case $h=\id_Y$.

\medskip

\begin{rem}
Let us note that if the endomorphism $h$ in \cref{th:comm_RB-end} is surjective, then $\mathfrak{R}\left(E\left(S\right) \right)= E\left(S\right)$. Indeed, if $\beta = h\left(\alpha\right)$ for a certain $\alpha\in Y$, observing that $\alpha\leq \beta$, we have that
$$
0_{\beta}
= 0_{h\left(\alpha\right)}
= \phii{h\left(\beta\right)}{h\left(\alpha\right)}\left(0_{h\left(\beta\right)} \right)
\underset{\eqref{eq:Worawiset}}{=} f_{\alpha}\phii{\beta}{\alpha}\left(0_{\beta}\right)
= f_{\alpha}\left(0_{\alpha}\right)
= \mathfrak{R}\left(0_{\alpha}\right),
$$
hence, by \cref{prop:RB-prop-E}, the claim follows.
\end{rem}

\bigskip

To conclude this section, we give a description of all idempotent Rota--Baxter endomorphisms of not necessarily abelian groups $\left(G,+\right)$. 
In this case, $G$ is exactly factorizable in $\ker\mathfrak{R}$ and $\im\mathfrak{R}$. 

\begin{theor}\label{prop_gruppi_RB}
   Let $(G,+)$ be a group. If $N\unlhd G$ is such that $G/N$ is abelian and $\mathcal{S}$ is a set of representatives of $G/N$ that is a subgroup of $G$, then any map $\mathfrak{R}:G\to G$ such that $\im\mathfrak{R}=\mathcal{S}$ and $\mathfrak{R}\left(g\right)\in N + g$, for every $g\in G$, is an idempotent Rota--Baxter endomorphism of $G$.\\
   Conversely, if $\mathfrak{R}$ is an idempotent Rota--Baxter endomorphism of $G$, then $G/\ker\mathfrak{R}$ is abelian, $\im\mathfrak{R}$ is a subgroup of $G$ that is also a set of representatives of $G/\ker\mathfrak{R}$, and $\mathfrak{R}\left(g\right)\in \ker\mathfrak{R} + g$, for every $g\in G$.
   \begin{proof}
   Initially, assume that $N\unlhd G$ and that $G/N$ is abelian, and let $\mathfrak{R}:G\to G$ be a map such that $\im\mathfrak{R} = \mathcal{S}$ and $\mathfrak{R}\left(g\right)\in N + g$, for every $g\in G$. Then, we clearly have that $\mathfrak{R}$ is a homomorphism from $G$ to $\mathcal{S}$. Besides, if $g,h\in G$, by \cite[Theorem 2.23]{Ro95}, since $G'\subseteq N$, 
   \begin{align*}
       \mathfrak{R}\left(g + \mathfrak{R}\left(g\right) + h - \mathfrak{R}\left(g\right)\right)
       &\in N + g + \mathfrak{R}\left(g\right) + h - \mathfrak{R}\left(g\right)
       = N + g + g + h - g\\
       &= N + g + [g,h] + h\subseteq N + g + h.
   \end{align*}
   Moreover, since $\mathfrak{R}\left(g\right) = n + g$, for a certain $n\in N$, we get
   $\mathfrak{R}^2\left(g\right) = \mathfrak{R}\left(n\right) + \mathfrak{R}\left(g\right)
       = \mathfrak{R}\left(g\right)$,
   hence the first part of the claim follows.\\
   Conversely, assume that $\mathfrak{R}:G\to G$ is an idempotent Rota--Baxter endomorphism of $G$. Then,  obviously, $\im \mathfrak{R}$ is a subgroup of $G$. Moreover, if $g,h\in G$, it holds that 
   \begin{align*}
      \mathfrak{R}\left(g + h - g - h\right)
      &= \mathfrak{R}\left(g\right) 
      - \mathfrak{R}\left(h\right) + \mathfrak{R}\left(h\right)  + \mathfrak{R}^2\left(h\right) + \mathfrak{R}\left(-g\right) - \mathfrak{R}^2\left(h\right)\\
      &= \mathfrak{R}\left(g\right) 
      - \mathfrak{R}\left(h\right)
      + \mathfrak{R}\left(h +\mathfrak{R}\left(h\right) - g - \mathfrak{R}\left(h\right)\right)\\
      &= \mathfrak{R}\left(g\right) - \mathfrak{R}\left(h\right) + \mathfrak{R}\left(h\right)  - \mathfrak{R}\left(g\right) = 0, 
   \end{align*}
   hence $g + h - g - h\in\ker\mathfrak{R}$, i.e., $G/\ker\mathfrak{R}$ is abelian.
   Now, if  $g\in G$, since
   \begin{align*}
      \mathfrak{R}\left(\mathfrak{R}\left(g\right)- g\right)
      = \mathfrak{R}^2\left(g\right) - \mathfrak{R}\left(g\right) = 0,
   \end{align*}
   we obtain that $\mathfrak{R}\left(g\right) = \mathfrak{R}\left(g\right) - g + g \in \ker\mathfrak{R} + g$. Thus, $\mathfrak{R}\left(g\right)\in \mathfrak{R}\left(G\right)\cap\left(\ker\mathfrak{R} + g\right)$. In addition, 
   if $\mathfrak{R}\left(h\right)\in\mathfrak{R}\left(G\right)\cap\left(\ker\mathfrak{R}+g\right)$, it follows that $\mathfrak{R}\left(g - h\right) = \mathfrak{R}\left(g\right) - \mathfrak{R}\left(h\right)\in\ker\mathfrak{R}$, hence
   $\mathfrak{R}\left(g - h\right) = \mathfrak{R}^2\left(g - h\right) = 0$, namely $g - h\in \ker\mathfrak{R}$. It follows that $g = n + h$, with $n\in \ker\mathfrak{R}$, thus
   \begin{align*}
       \mathfrak{R}\left(g\right)
       = \mathfrak{R}\left(n\right) + \mathfrak{R}\left(h\right)
       = \mathfrak{R}\left(h\right).
   \end{align*}
   Therefore, $\im\mathfrak{R}$ is a set of representatives of $G/\ker \mathfrak{R}$. 
   \end{proof}
\end{theor}

\bigskip

\section{Dual weak braces obtained by Rota--Baxter operators}
\noindent In this section, we illustrate a method for obtaining a dual weak brace $S$ starting from a given Rota--Baxter operator on a Clifford semigroup $(S,+)$. This construction is inspired by that of skew braces due to Bardakov and Gubarev \cite[Proposition 3.1]{BaGu21}. Furthermore, we show that every commutative Rota--Baxter endomorphism determines a bi-weak brace.

\medskip

\begin{lemma}
\label{th:prod-Clifford}
  Let $(S,+)$ be a Clifford semigroup and $\mathfrak{R}: S \to S$ a Rota--Baxter operator on $(S,+)$. Set
    $a\circ b:= a + \mathfrak{R}\left(a\right) + b - \mathfrak{R}\left(a\right)$,
    for all $a, b \in S$. Then, $(S, \circ)$ is a Clifford semigroup.
    \begin{proof}
  It is a routine computation  to verify that  $(S, \circ)$ is a semigroup. Moreover, if $a\in S$, set $x:=-\mathfrak{R}\left(a\right)-a+\mathfrak{R}\left(a\right)$, we get
\begin{align*}
    a \,\circ  x\circ a
    &=\left(a+\mathfrak{R}\left(a\right)^0-a+\mathfrak{R}\left(a\right)^0\right)\circ a\\
    &= a^0\circ a &\mbox{by \cref{prop:RB-prop-E}}\\
    &=a^0 +a+ \mathfrak{R}\left(a^0\right)-\mathfrak{R}\left(a^0\right)&\mbox{by \cref{prop:RB-prop-E}}\\
    &= a^0+a =a.
\end{align*}
To prove the claim we use \cite[Theorem 4.2.1]{Ho95} and thus
we show that the idempotents of $(S, \circ)$ are central.
To this end, observe that $E(S, \circ)\subseteq E(S, +)$. 
Indeed, if $e \in E(S, \circ)$, then
\begin{align*}
  \mathfrak{R}\left(e\right) + \mathfrak{R}\left(e\right)
  = \mathfrak{R}\left(e + \mathfrak{R}\left(e\right) + e - \mathfrak{R}\left(e\right)\right)
  = \mathfrak{R}\left(e\circ e\right)
  = \mathfrak{R}\left(e\right),
\end{align*}
i.e., $\mathfrak{R}\left(e\right)\in E(S, +)$, and thus
\begin{align*}
    e = e\circ e = e +\mathfrak{R}\left(e\right)+e-\mathfrak{R}\left(e\right)
    = e + e +\mathfrak{R}\left(e\right)^0
    = e + e.
\end{align*}
Hence, if $a \in S$, we get
\begin{align*}
    a \circ e&=a+\mathfrak{R}\left(a\right)+e-\mathfrak{R}\left(a\right)=a+\mathfrak{R}\left(a\right)^0+e=a+e\\
    &= a + e+ \mathfrak{R}\left(e\right)^0
    = e + \mathfrak{R}\left(e\right) + a -\mathfrak{R}\left(e\right)
    = e\circ a.
\end{align*}
Therefore, $(S, \circ)$ is a Clifford semigroup and $a^-$ coincides with $x$.
\end{proof}
\end{lemma}   

\medskip

\begin{rem}
  As it happens in the context of groups (see \cite[Proposition 3.1]{GuLaSh20}), if $\mathfrak{R}$ is a Rota--Baxter operator on a Clifford semigroup  $(S,+)$, then clearly 
$\mathfrak{R}$ is a homomorphism from $(S, \circ)$ to $(S,+)$. Moreover, $\mathfrak{R}$ is also a Rota--Baxter operator on $(S, \circ)$, indeed if $a,b\in S$ we have that
  \begin{align*}
      \mathfrak{R}\left(a\right)\circ\mathfrak{R}\left(b\right) 
          &= \mathfrak{R}\left(a\right) + \mathfrak{R}^2\left(a\right)
          + \mathfrak{R}\left(b\right)
          - \mathfrak{R}^2\left(a\right)
          = \mathfrak{R}\left(a\circ\mathfrak{R}\left(a\right)
          \circ b\right)
          + \mathfrak{R}\left(\mathfrak{R}\left(a\right)^-\right)\\
          &= \mathfrak{R}\left(a\circ\mathfrak{R}\left(a\right)
          \circ b\circ\mathfrak{R}\left(a\right)^-\right)
  \end{align*}
  and, by \eqref{eq:circ+},
  \begin{align*}
      a\circ \mathfrak{R}\left(a\right)\circ  \mathfrak{R}\left(a\right)^-
      &= a +\lambda_a\mathfrak{R}\left(a\right)-\lambda_a\mathfrak{R}\left(a\right)
      = a\circ \mathfrak{R}\left(a\right)
      - a\circ \mathfrak{R}\left(a\right) + a\\
      &= a^0\circ\mathfrak{R}\left(a\right)^0 + a
      =  a + a^0 + \mathfrak{R}\left(a\right)^0 = a.
  \end{align*}
  \medskip
\end{rem}

\medskip

\begin{theor}\label{teo_skewsemibrace}
   Let $(S,+)$ be a Clifford semigroup. If $\mathfrak{R}: S \to S$ is a Rota--Baxter operator on $(S,+)$,  then $(S,+)$ endowed with the operation $\circ$  defined by
       \begin{align}\label{def: cerchietto}
        a\circ b:= a + \mathfrak{R}\left(a\right) + b - \mathfrak{R}\left(a\right),
    \end{align}
    for all $a, b \in S$, is a dual weak brace, which we call \emph{weak brace associated to $\mathfrak{R}$} and denote by $S_\mathfrak{R}$. 
\begin{proof}
Initially, by \cref{th:prod-Clifford}, $\left(S,\circ\right)$ is a Clifford semigroup where, for any $a\in S$, $a^- = -\mathfrak{R}\left(a\right)-a+\mathfrak{R}\left(a\right)$.
Moreover, if $a,b,c\in S$, we get
\begin{align*}
    a\circ b - a + a\circ c
    &= a + \mathfrak{R}\left(a\right) + b - \mathfrak{R}\left(a\right) + a^0 + \mathfrak{R}\left(a\right) + c - \mathfrak{R}\left(a\right)\\
    &= a + \mathfrak{R}\left(a\right) + b  +\mathfrak{R}\left(a\right)^0+ c - \mathfrak{R}\left(a\right)\\
    &= a + \mathfrak{R}\left(a\right) + b + c - \mathfrak{R}\left(a\right) 
    = a\circ\left(b + c\right),
\end{align*}
In addition, it holds $a \circ a^-
    = a + \mathfrak{R}(a)^0-a+\mathfrak{R}\left(a\right)^0
    = a^0$.
Therefore, the statement is proved.
\end{proof}
\end{theor}

\medskip

\noindent For example, the dual weak braces associated to the elementary Rota--Baxter operators $\mathfrak{E}$ and $\mathfrak{O}$ on a Clifford semigroup $(S,+)$ are the trivial and the almost trivial weak braces, respectively. Moreover, if $(S,+)$ is commutative, then any weak brace $S_{\mathfrak{R}}$ related to an arbitrary Rota--Baxter operator $\mathfrak{R}$ on $S$ is a trivial weak brace. 

\medskip

Let us note that if  $\mathfrak{R}$ and $\mathfrak{T}$ are two equivalent Rota--Baxter operators on two Clifford semigroups $(S,+)$ and $(T,+)$, respectively, then, $S_{\mathfrak{R}}$ and $T_{\mathfrak{T}}$ are isomorphic. Indeed, if $\varphi$ is the isomorphism from $(S,+)$ to $(T,+)$ such that $\mathfrak{R}=\varphi^{-1}\mathfrak{T}\varphi$, then
\begin{align*}
\varphi(a\circ b) &=  \varphi(a+\mathfrak{R}(a)+b-\mathfrak{R}(a)) = \varphi(a) + \varphi\mathfrak{R}(a)+\varphi(b)-\varphi\mathfrak{R}(a)\\
&= \varphi(a) + \mathfrak{T}\varphi(a)+\varphi(b)-\mathfrak{T}\varphi(a)
= \varphi(a) \circ \varphi(b),
\end{align*}
for all $a,b\in S$.\\ 
It is natural to ask whether the converse holds. The answer is negative, since if we consider a commutative Clifford semigroup $(S,+)$, with $E(S) \neq S$, the trivial weak brace on $(S,+)$ can be obtained starting from arbitrary Rota--Baxter operators on $(S,+)$.\\
Analogously to \cite[Question 7.4]{BaGu21}, the following question arises.
\begin{que}
Given two Rota--Baxter operators $\mathfrak{R}$ and $\mathfrak{T}$ on a Clifford semigroup $(S,+)$, under which conditions the associated weak braces $S_\mathfrak{R}$ and $S_\mathfrak{T}$ are isomorphic?
\end{que}

\medskip

 If $(S,+)=[Y; G_{\alpha}; \phi_{\alpha, \beta}]$ is a Clifford semigroup in which each group $G_{\alpha}$ has a Rota--Baxter operator $\mathfrak{R}_{\alpha}$, the strong semilattice of the skew braces $G_{\mathfrak{R}_{\alpha}}$ coincides with the dual weak brace $S_{\mathfrak{R}}$, where $\mathfrak{R}$ is the strong Rota--Baxter operator on $(S,+)$, as we show below.

\begin{rem}
Let $Y$ be a (lower) semilattice, $\{G_{\alpha} \ | \ \alpha\in Y\}$, a family of disjoint groups each one having a Rota--Baxter operator $\mathfrak{R}_{\alpha}$. 
Thus, by \cref{teo_skewsemibrace}, $\{G_{\alpha} \ | \ \alpha\in Y\}$ is a family of disjoint skew braces associated to $\mathfrak{R}_{\alpha}$, for every $\alpha\in Y$. Moreover, considered a strong semilattice $S=[ Y; G_{\alpha}; \phii{\alpha}{\beta}]$ of the groups $G_{\alpha}$,  by \cref{prop:Strong-Lattice-Inverse-Semi-Brace}, the structure $\left(S,+, \circ\right)$ is a dual weak brace, with $a+b= \phii{\alpha}{\alpha\beta}(a)+\phii{\beta}{\alpha\beta}(b)$ and
\begin{align*}
     a\circ b= \phii{\alpha}{\alpha\beta}(a)\circ\phii{\beta}{\alpha\beta}(b)=\phii{\alpha}{\alpha\beta}\left(a\right)+\mathfrak{R}_{\alpha\beta}\phii{\alpha}{ \alpha\beta}\left(a\right)+\phii{\beta}{ \alpha\beta}\left(b\right)-\mathfrak{R}_{\alpha\beta}\phii{\alpha}{ \alpha\beta}\left(a\right),
\end{align*}
for all $a\in G_{\alpha}$ and $b\in G_{\beta}$.
On the other hand, by \cref{teo_RB_semilattice}, we have that the map $\mathfrak{R}:S \to S$ given by $\mathfrak{R}\left(a\right) = \mathfrak{R}_{\alpha}\left(a\right)$, for every $a\in G_{\alpha}$, is a Rota--Baxter operator on the Clifford semigroup $\left(S, +\right)$ if and only if \eqref{condizione_semireticolo} is satisfied. Hence, in this case, by \cref{teo_skewsemibrace}, $S_{\mathfrak{R}}=(S, +, \circ_{\mathfrak{R}})$ is a weak brace and
\begin{align*}
  a \circ_{\mathfrak{R}} b
  &= a + \mathfrak{R}_{\alpha}\left(a\right) + b - \mathfrak{R}_{\alpha}\left(a\right)\\
  &= \phii{\alpha}{\alpha\beta}\left(a\right) + \phii{\alpha}{\alpha\beta}\mathfrak{R}_{\alpha}\left(a\right) + \phii{\beta}{\alpha\beta}\left(b\right) - \phii{\alpha}{\alpha\beta}\mathfrak{R}_{\alpha}\left(a\right)\\
  &= \phii{\alpha}{\alpha\beta}\left(a\right)+\mathfrak{R}_{\alpha\beta}\phii{\alpha}{ \alpha\beta}\left(a\right)+\phii{\beta}{ \alpha\beta}\left(b\right)-\mathfrak{R}_{\alpha\beta}\phii{\alpha}{ \alpha\beta}\left(a\right) &\mbox{by \eqref{condizione_semireticolo}}\\
  &
  = a\circ b,
\end{align*}
for all $a\in G_{\alpha}$ and $b\in G_{\beta}$.
\end{rem}

\medskip

Instances of bi-weak braces can be obtained through special Rota--Baxter operators, as we show in the following result that is related to \cite[Proposition 5.1]{Ko21}.

\begin{prop}
Let $(S,+)$ be a Clifford semigroup and consider a commutative Rota--Baxter endomorphism $\mathfrak{R}$ of $(S,+)$. Then, the weak brace $S_{\mathfrak{R}}$ is a bi-weak brace.
\begin{proof}
If $a,b,c \in S$. Since $\mathfrak{R}\left(a^-\right)=-\mathfrak{R}\left(a\right)$, we get
\begin{align*}
    \left(a + b\right) \circ a^- \circ \left(a +c\right)
    &=\left(a+b\right) \circ \left(-\mathfrak{R}\left(a\right)+a^0+c+\mathfrak{R}\left(a\right)\right)\\
&=a+b+\mathfrak{R}\left(a\right)+\mathfrak{R}\left(b\right)-\mathfrak{R}\left(a\right)+a^0+c+\mathfrak{R}\left(a\right)-\mathfrak{R}\left(b\right)-\mathfrak{R}\left(a\right)\\
&=a+b+\mathfrak{R}\left(a\right)^0+\mathfrak{R}\left(b\right)+c-\mathfrak{R}\left(b\right)+\mathfrak{R}\left(a\right)^0\\
&=a+b+\mathfrak{R}\left(b\right)+c-\mathfrak{R}\left(b\right)\\
&=a + ( b \circ c).
\end{align*}
Moreover, since $S_{\mathfrak{R}}$ is a dual weak brace, it trivially holds $a^- \circ a=a-a$. Therefore, the claim follows.
\end{proof}
\end{prop}

\bigskip

\section{Solutions associated to Rota--Baxter operators}
In this section, we make explicit the solutions associated to dual weak braces obtained through Rota-Baxter operators. In particular, such maps have a behaviour near to bijectivity and non-degeneracy.

\medskip

Let us recall that every weak brace $(S,+, \circ)$ gives rise to a solution $r:S\times S\to S\times S$ defined by
\begin{align*}
    r\left(a,b\right)
    = \left(-a + a\circ b, \ \left(-a + a\circ b\right)^-\circ a\circ b\right),
\end{align*}
for all $a,b\in S$ (see \cite[Theorem 11]{CaMaMiSt21x}).
In particular, the images of $r$ can be written by the maps $\lambda_a$ and $\rho_b:S\to S$ given by 
\begin{align*}
    \rho_b\left(a\right)= \lambda_a\left(b\right)^-\circ a \circ b,
\end{align*}
for all $a,b\in S$.  The map $\rho:S\to \Map(S), b\mapsto \rho_b$  is a semigroup anti-homomorphism of the inverse semigroup $\left(S,\circ\right)$ into the monoid $\Map(S)$ of the maps from $S$ into itself.\\
Any solution $r$ associated
to an arbitrary weak brace $(S, +, \circ)$ has a behaviour close to bijectivity, namely $r$ is a
completely regular element in $\Map(S \times S)$. Indeed, considered the solution $r^{op}$ associated to the \emph{opposite weak brace}
$S^{op}$ of $S$, i.e., the structure $\left(S,+^{op}, \circ\right)$ where $a+^{op}b:= b + a$, for all $a,b \in S$, they hold
  \begin{align*}
      r\, r^{op}\, r = r, \qquad
      r^{op}\, r\, r^{op} = r^{op}, \qquad \text{and}\qquad rr^{op} = r^{op}r.
  \end{align*}
In general, the maps $\lambda_a$ and $\rho_b$ are not completely regular in
$\Map(S)$ (see \cite[Example 3]{CaMaMiSt21x}), unless $S$ is a dual weak brace. Indeed, in this case, $r$ has also a behaviour close to the non-degeneracy because they hold
\begin{align*}
    \lambda_a\lambda_{a^-}\lambda_a
    = \lambda_{a},
    \qquad   \lambda_{a^-}\lambda_{a}\lambda_{a^-}
    = \lambda_{a^-}, \quad &\text{and} \qquad
    \lambda_a\lambda_{a^-}
    =\lambda_{a^-}\lambda_{a},\\
    \rho_{b}\rho_{b^-}\rho_{b}
    = \rho_b,
    \qquad
    \rho_{b^-}\rho_{b}\rho_{b^-} = \rho_{b^-},
    \qquad &\text{and} \qquad
    \rho_{b}\rho_{b^-}
    = \rho_{b^-}\rho_{b},
\end{align*}
for all $a, b \in S$.
In particular, if $S$ is a skew brace, $r$ is bijective with  $r^{-1} = r^{op}$ and $r$ is non-degenerate.

\medskip

Note that if $(S,+, \circ)$ is a dual weak brace and the map $\lambda_a$ is bijective, for every $a\in S$, then $S$ is a skew brace. Hence, since the maps $\lambda_a$ are completely regular, we obtain that $\lambda_{a}^{-1}=\lambda_{a^-}$, for every $a\in S$. Thus, $\lambda_e=\id_S$, for every $e \in E(S)$, and if $e_1,e_2\in E(S)$, we get
\begin{align*}
    e_1=\lambda_{e_2}(e_1)=-e_2+e_2\circ e_1=e_2+e_2+e_1=e_2+e_1.
\end{align*}
Similarly, $e_2=e_1+e_2$, and so $e_1=e_2$. Therefore, all the idempotents are equal and $|E(S)|=1$.

\bigskip

Now, let us make explicit the solution $r$ associated to a dual weak brace $S_{\mathfrak{R}}$ associated to a Rota--Baxter operator on a Clifford semigroup $(S,+)$, that is 
\begin{align}\label{cor:sol-SRB}
    r\left(a,b\right)
    = \left(\lambda_a\left(b\right), 
    \ -\mathfrak{R}\lambda_a\left(b\right)-\lambda_a\left(b\right) + a + \lambda_a\left(b\right) +\mathfrak{R}\lambda_a\left(b\right)\right),
\end{align}
for all $a,b\in S$. In fact,
\begin{align*}
    \rho_b\left(a\right)&=\lambda_a\left(b\right)^-\circ \left(a \circ b\right)=\lambda_a\left(b\right)^-+\mathfrak{R}\left(\lambda_a\left(b\right)^-\right)+a \circ b-\mathfrak{R}\left(\lambda_a\left(b\right)^-\right)\\
    &=\lambda_a\left(b\right)^--\mathfrak{R}\lambda_a\left(b\right)+a+\lambda_a(b)+\mathfrak{R}\lambda_a\left(b\right) &\mbox{by \eqref{eq:circ+}}\\
    &=-\mathfrak{R}\lambda_a\left(b\right)-\lambda_a(b)+\mathfrak{R}\lambda_a\left(b\right)^0+a+\lambda_a(b)+\mathfrak{R}\lambda_a\left(b\right)\\
    &=-\mathfrak{R}\lambda_a\left(b\right)-\lambda_a\left(b\right) + a + \lambda_a\left(b\right) +\mathfrak{R}\lambda_a\left(b\right),
\end{align*}
for all $a,b \in S$.

\medskip

Now, given a Rota--Baxter $\mathfrak{R}$ on a Clifford semigroup $(S,+)$, we show that the solution $r^{op}$ associated to the opposite weak brace $S_\mathfrak{R}^{op}$  of $S_\mathfrak{R}$ is strictly linked to the solution associated to $S_{\mathfrak{R}^{op}}$, where $\mathfrak{R}^{op}$ is the opposite Rota--Baxter operator of $\mathfrak{R}$ defined by $\mathfrak{R}^{op}(a)=-a+\mathfrak{R}(-a)$, for every $a \in S$. To this purpose, we recall that two solutions $r,s$ on a set $S$ are said to be \emph{equivalent} (or \emph{isomorphic}) if there exists a bijection $f$ on $S$ such that $\left(f\times f\right)r = s\left(f\times f\right)$, see \cite[p. 105]{CeJeOk14}.

\begin{prop}
Let $S_\mathfrak{R}$ be the dual weak brace associated to a Rota--Baxter operator $\mathfrak{R}$ on a Clifford semigroup $\left(S,+\right)$ and $r$ the solution associated to $S_\mathfrak{R}$.
Then,
 $S_{\mathfrak{R}^{op}}=(S, +, \tilde{\circ})$ and  $S_{\mathfrak{R}}^{op}$ are isomorphic and
   the opposite solution $r^{op}$ of $r$ is equivalent to the solution associated to $S_{\mathfrak{R}^{op}}$. In particular, it  is given by
    \begin{align*}
        r^{op}(a,b)&=\left(\lambda_a\left(b\right)^{op}, \ -\mathfrak{R}\left(\lambda_a\left(b\right)^{op}\right)-\lambda_a\left(b\right)^{op}+a+\lambda_a\left(b\right)^{op}+\mathfrak{R}\left(\lambda_a\left(b\right)^{op}\right)
        \right),
    \end{align*}
    where $\lambda_a(b)^{op}=a + \mathfrak{R}\left(a\right) + b -\mathfrak{R}\left(a\right) - a$, for all $a,b \in S$.
\begin{proof}
Let $a,b\in S$, then
\begin{align*}
    a\ \tilde{\circ} \ b 
    =a^0+ \mathfrak{R}\left(-a\right)
    + b - \mathfrak{R}\left(-a\right) + a= \mathfrak{R}\left(-a\right)
    + b - \mathfrak{R}\left(-a\right) + a
    = -\left(\left(-a\right)\circ \left(-b\right)\right).
\end{align*}
Besides, on can check that $f:S\to S$ defined by $f\left(a\right) = -a$, for every $a\in S$, is an isomomorphism of weak braces from $S_{\mathfrak{R}^{op}}$ to $S^{op}_{\mathfrak{R}}$.\\
Now, it is a routine computation to calculate the solution $r^{op}$ as in \cite[Remark 2]{CaMaMiSt21x} and, if $\tilde{r}$ is the solution associated to $S^{op}_{\mathfrak{R}}$ as in \eqref{cor:sol-SRB}, then $f$ is a bijection of $S$ satisfying
$\left(f\times f\right)r^{op} = \tilde{r}\left(f\times f\right)$, namely,  $r^{op}$ and $\tilde{r}$ are equivalent.
\end{proof}
\end{prop}

\bigskip

\section{Ideals of dual weak braces}

\noindent In this section, we deal with the ideals of a dual weak brace, by extending the already known theory for skew braces. In particular, we introduce the notion of the socle of a dual weak brace. 
Finally, we describe the ideals of dual weak braces associated to a given Rota--Baxter operator.
\medskip

Inverse semigroups are a widely studied class of semigroups for which congruences are well understood. Indeed, in this context it is possible to mimic the group-theoretical treatment of congruences since one can find special subsemigroups that are analogue of normal subgroups of groups. 
To this end, we recall the definition of such structures below.
\begin{defin}\label{normal}\cite[Definition VI.1.2]{PeRe99} 
Let $S$ be a Clifford semigroup. A subset $N$ of $S$ is  \emph{normal} if it satisfies the following conditions: 
    \begin{enumerate}
        \item $E(S)\subseteq N$,
        \item $\forall \ a\in S\quad a\in N\Longrightarrow a^{-1}\in N$, 
        \item $\forall \ a,b\in S\quad a,\; a^{0}b\in N\Longrightarrow ab\in N$, 
        \item $\forall \ a,b\in S\quad ab\in N\Longrightarrow ba\in N$.
    \end{enumerate}
\end{defin}
\noindent Note that a subset $N$ of a group $G$ is a normal subgroup of $G$ if and only if it is a normal subset of $G$ (see \cite[Exercises VI.1.22(i)]{PeRe99}. Given a Clifford semigroup $S=[Y; G_{\alpha}; \phii{\alpha}{\beta}]$, for each $\alpha\in Y$, let $K_\alpha$ be a normal subgroup of $G_\alpha$ and assume that $\phii{\alpha}{\beta}(K_\alpha)\subseteq K_\beta$ if $\alpha > \beta$. Set $\psi_{\alpha,\beta} = {\phi_{\alpha,\beta}}_{|_{K_\alpha}}$ if $\alpha \geq \beta$ and $K=[Y; K_{\alpha}; \psi_{\alpha, \beta}]$, then $K$ is a normal subset of $S$; conversely, every normal subset of $S$ is of this form (see \cite[Exercises III.1.9(ii)]{Pe84}).
\\
Clearly, if $S$ is a Clifford semigroup, then $E(S)$ is a normal subset of $S$. Indeed, the conditions (1), (2), and (3) are trivially satisfied; moreover, the condition (4) follows by \cite[Lemma II.2.2(iii)]{PeRe99}.
\medskip

\begin{rems}\label{rem:norm} 
Let $N$ be a normal subset of a  Clifford semigroup $S$. 
\begin{enumerate}
    \item Observe that $N$ is a subsemigroup of $S$. Indeed, if $a,b \in N$, then, by $1.$ of \cref{normal}, we get $b^{0}a^0\in N$ and, by $3.$, $ba^0 \in N$, thus $a^0b\in N$. Hence, by $3.$, it follows that $ab\in N$.
    \vspace{1mm}
    \item If $n\in N$ and $a\in S$, then $a^{-1}na\in N$. In fact, since $a^0\in N$, by $1.$, we have that $\left(na\right)a^{-1}\in N$ and, by $4.$ in \cref{normal}, $a^{-1} n a\in N$.
\end{enumerate}
\end{rems}

\medskip

The following definition is essential for obtaining a quotient structure and it is consistent with that given for skew braces (see \cite[Definition 2.1]{GuVe17}).
\begin{defin}\label{def:ideal}
Let $(S,+, \circ)$ be a dual weak brace. Then, a subset $I$ of $S$ is an \emph{ideal} of $S$ if the following hold:
    \begin{enumerate}
        \item $I$ is a normal subset of $\left(S, +\right)$,
        \item $\lambda_a\left(I\right)\subseteq I$, for every $a\in S$, 
        \item $I$ is a normal subset of $\left(S, \circ\right)$.
    \end{enumerate}
\end{defin}

\medskip

\noindent Observe that if $S$ is a dual weak brace, then $S$ and $E(S)$ are trivial ideals of $S$.

\begin{theor}\label{teo_congru_skewinverse}
Let $I$ be an ideal of a dual weak brace $(S,+,\circ)$. Then, the relation $\sim_{I}$ on $S$ given by
\begin{align*}
    \forall \ a,b\in S\quad a\sim_I b  
    \ \iff \
    a^0 = b^0
     \ \  \text{and} \ \ -a+ b\in I,
\end{align*}
is a congruence of $(S,+, \circ)$. Moreover,  $S/\sim_I$ is a dual weak brace with semilattice of idempotents isomorphic to $E\left(S\right)$.
\begin{proof}
Initially, we prove that 
\begin{align*}
     a\sim_I b  
    \ \iff \
    a^0 = b^0
     \ \  \text{and} \ \ a^-\circ b\in I,
\end{align*}
for all $a,b \in S$. Indeed, if $a,b\in S$ are such that $a\sim_I b$, then, clearly $a^0=b^0$, and, since $\lambda_x\left(x^-\right)=-x$, for every $x \in S$ (see \cite[Lemma 2]{CaMaMiSt21x}), we get
$$
   a^- \circ b=a^-+\lambda_{a^-}\left( b\right)=\lambda_{a^-}\left(-a+b\right) \in I.$$
Conversely, if $a,b \in S$ are such that $a^0=b^0$ and $a^-\circ b\in I$, then
\begin{align*} 
-a+b=\lambda_a\left(a^-\right)+\lambda_{b^0}\left(b\right)=\lambda_a\left(a^-\right)+\lambda_{a^0}\left(b\right)=\lambda_a\left(a^-+\lambda_{a^-}\left(b\right)\right)
    = \lambda_a\left(a^-\circ b\right) \in I.
\end{align*}
Now, by \cite[Lemma VI.3.1]{PeRe99}, \cite[Exercise VI.2.13(i)]{PeRe99}, and \cite[Theorem 4.2.1]{Ho95}, $\sim_I$ is a congruence of the dual weak brace $\left(S, +, \circ \right)$ whose both quotient structures are Clifford semigroups. The remaining parts of the statement are for immediate verification, so the proof is complete.
\end{proof}
\end{theor}

\medskip

The following result describes the ideals of dual weak braces associated to Rota--Baxter operators and it consistent with  \cite[Proposition 5.4(a)]{BaGu21}. Into the specific, the request of $\lambda$-invariance is redundant.
\begin{prop}
Let $S_{\mathfrak{R}}$ be the dual weak brace obtained by a Rota--Baxter operator $\mathfrak{R}$ on a Clifford semigroup $(S,+)$ and $I\subseteq S$. Then, $I$ is an ideal of $S_{\mathfrak{R}}$ if and only if $I$ both is a normal subset of the Clifford semigroups $(S,+)$ and $(S,\circ)$.
\begin{proof}
Assume that $I$ is a normal subset both of $(S,+)$ and $(S,\circ)$. Thus, $-a+x+a\in I$, for all $x\in I$ and $a\in S$. Indeed, since $a^0\in I$, it follows that $a^0+x=a +(-a + x)\in I$ and so, by $4.$ in \cref{normal}, $-a + x + a\in I$. 
Hence, if $a\in S$ and $x\in I$, we obtain that
\begin{align*}
    \lambda_a\left(x\right)
    = -a + a\circ x
    = a^0 + \mathfrak{R}\left(a\right) + x -\mathfrak{R}\left(a\right)\in I.
\end{align*}
Therefore, since the converse is trivial, the claim follows.
\end{proof}
\end{prop}

\bigskip

If $S$ is a dual weak brace, we call \emph{socle} of $S$ the following set
 \begin{align*}
    \Soc\left(S\right)=\{ a  \, | \,a \in S, \, \,\forall \ b \in S \quad a+b=a \circ b \quad \text{and} \quad a+b=b+a \}.
\end{align*}  
 
\noindent Clearly, the previous definition includes that already known in the context of skew braces in \cite{GuVe17} originally introduced in \cite{Ru05}. 
Moreover, by the following result, we can always consider the quotient structure of a dual weak brace by its socle.
\begin{prop}
Let $(S,+, \circ)$ be a dual weak brace. Then, $\Soc(S)$ is an ideal of $S$.
\begin{proof}
Initially, note that $E(S) \subseteq \Soc(S)$. Indeed, if $e \in E(S)$ and $a \in S$, it is sufficient to show that
\begin{align}\label{idemp}
    e+a=e \circ a.
\end{align}
In fact, $
    e \circ a=e +\lambda_e\left(a\right)=e-e+e\circ a=\lambda_e\left(a\right),
$
and, by \cref{+circ}, 
\begin{align*}
e+a=e \circ \lambda_{e^-}\left(a\right)=e \circ e^- \circ a=e \circ a.
\end{align*}
Now, let us prove that $\lambda_x\left(\Soc\left(S\right)\right) \subseteq \Soc\left(S\right)$, for every $x \in S$. If $a \in \Soc(S)$ and $b \in S$, we obtain
\begin{align*}
    \lambda_x(a) \circ b&=x \circ \left(a+x^-\right) \circ b=x \circ \left( a \circ x^- \circ b\right)=x \circ \left( a + x^- \circ b\right)\\
    &=x \circ a -x +x^0 \circ b\\
    &=x \circ a -x +x^0 + b&\mbox{by \eqref{idemp}}\\
    &= x \circ (a+x^-) +b=\lambda_x(a)+ b
\end{align*}
and, on the other hand, by the previous computations,
\begin{align*}
    \lambda_x(a)+ b&=\lambda_x(a) \circ b=x \circ  a \circ x^- \circ b=b^0 \circ x \circ \left( a \circ x^- \circ b\right)\\
    &=b^0 \circ x \circ  \left(a + x^- \circ b\right)=b \circ \left( b^- \circ x \circ  \left( x^- \circ b + a\right)\right)\\
    &=b \circ \lambda_{b^- \circ x}\left(a\right)\\
    &= b+ \lambda_x(a)& \mbox{by \cref{+circ}}
\end{align*}
Now, it is a routine computation to prove that $\Soc(S)$ is closed with respect to the $\circ$. Hence, if $a, a^0 \circ b \in \Soc(S)$, then $a \circ b=a \circ a^0 \circ b \in Soc(S)$. Moreover, if $a \in \Soc(S)$, then $a^- \in \Soc(S)$. Indeed, we have
\begin{align*}
    a^-+b&= a^- \circ \left(-a+a\circ b\right) & \mbox{by \cref{+circ}}\\
    &=a^- \circ \left(-a+a+ b\right)\\
    &=a^- \circ a^0 \circ b  &\mbox{by \eqref{idemp}}\\
    &=a^- \circ b.
\end{align*}
Considering $b=a$, it follows that 
\begin{align}\label{menoa}
    a^-+a=a^0.
\end{align}
Thus, we get
\begin{align*}
    b+a^-&=b+a^- \circ \left(a^-\right)^0=b+a^- \circ a^0\\
    &=b+a^- + a^0 &\mbox{by \eqref{idemp}}\\
    &=a^- \circ a +b+a^-\\
    &=a^-+a+b+a^- &\mbox{by \eqref{menoa}}\\
    &=a^-+b+a+a^-\\
    &=a^-+b+a^0&\mbox{by \eqref{menoa}}\\
    &=a^-+b,
\end{align*}
and so $a^- \in \Soc (S)$. It immediately follows that $-a=\lambda_a\left(a^-\right) \in \Soc(S)$. Finally, if $a\circ b \in \Soc(S)$,
\begin{align*}
    b \circ a=b \circ \left(a \circ b \circ b^-\right)=b \circ \left( a \circ b +b^-\right) =b \circ \left(b^-+ a \circ b \right) =\lambda_b(a \circ b),
\end{align*}
thus $b \circ a \in \Soc(S)$, i.e., $\Soc(S)$ is a normal subset of $(S, \circ)$.\\
It is easy to prove that $\Soc(S)$ is closed with respect to the $+$, hence if $a, a^0+b \in \Soc (S)$, clearly $a+b \in \Soc(S)$. To conclude the proof, we observe that if $a+b\in \Soc(S)$, then 
\begin{align*}
    b+a=-a+(a+b)+a=-a+a+(a+b)=a+b \in \Soc(S).
\end{align*}
Therefore, the statement is proved.
\end{proof}
\end{prop}

\medskip

The following proposition provides a method to obtain ideals of dual weak braces starting from two given ones.
\begin{prop}
    Let $S$ be a dual weak brace and $I, J$ ideals of $S$. Then, $I+J$ and $I\circ J$ are ideals of $S$.
    \begin{proof}
    Initially, by $2.$ in \cref{rem:norm}, one can check that $I+J=\{i + j \ | \ i\in I \ \text{and} \ j\in J\}$ is closed with respect to $+$. Moreover, $\lambda_a(i+j)=\lambda_a(i)+\lambda_a(j)\in I+J$, for every $a \in S$. Hence, by \eqref{eq:circ+}, it also follows that $I+J$ is closed with respect to $\circ$ since,  if $a,b\in I+J$,
    $a\circ b = a + \lambda_{a}\left(b\right)\in I+J$. In addition, it is clear that $E(S) \subseteq I+J$. \\
    Now, let us show that $I+J$ is a normal subset of $(S,+)$.  If $a = i + j\in I+J$, then
    $-a = \left(-j-i+j\right) - j\in I + J$.
    Moreover, if $a\in I+J$ and $a^0+b\in I+J$, then
    \begin{align*}
        a + b = a +  a^0 + b\in I+J
    \end{align*}
   and if $a+b = i+j\in I+J$, it follows that
   \begin{align*}
       b + a 
       = a^0 + b +a
       = -a + \left(i + j\right) + a
       =-a + i + a -a + j + a\in I + J.
   \end{align*}
  Thus, $I+J$ is a normal subset of $(S,+)$. Finally, we prove that $I+J$ is a normal subset of $(S,\circ)$. Now, observe that $a^-= -\lambda_{a^-}\left(a\right)\in I+J$ and, if $a, a^0\circ b\in I+J$, then
  \begin{align*}
      a\circ b = a\circ a^0\circ b\in I + J.
  \end{align*}
  Furthermore, if $a\circ b = i + j\in I+J$, by \eqref{eq:circ+} and \cref{+circ}, we obtain that
  \begin{align*}
     b \circ a
     = a^- \circ (i +j )\circ a 
     = a^-\circ i\circ\lambda_{i^-}\left(j\right)\circ a
     = \left(a^-\circ i\circ a\right)\circ
     \left(a^-\circ\lambda_{i^-}\left(j\right)\circ a\right)\in I+J,
  \end{align*}
 therefore $I+J$ is an ideal of $S$.\\
To prove the second part of the statement, note that $I\circ J=\{i\circ j \ | \ i\in I \ \text{and} \ j\in J\}$ is closed with respect to $\circ$. 
    Besides, if $a,b\in I\circ J$, it follows that
    $a + b = a\circ\lambda_{a^-}\left(b\right)\in I\circ J$, hence $I\circ J$ is closed with respect to $+$.
    Moreover,  by \eqref{eq:circ+}, we  have that $\lambda_a(i\circ j)=\lambda_a(i)+\lambda_{a\circ i}(j)\in I\circ J$, for every $a \in S$. 
    In addition, it is clear that $E(S) \subseteq I+J$. \\
    Now, let us show that $I\circ J$ is a normal subset of $(S,+)$.  If $a = i\circ j\in I\circ J$, then
     $-a = \lambda_{i}\left(-j\right) - i 
     = -i + \left(i+\lambda_{i}\left(-j\right) - i\right)\in I\circ J$.
     Moreover, if $a\in I\circ J$ and $a^0+b\in I\circ J$, then
   $
        a + b = a +  a^0 + b\in I\circ J
  $
   and if $a+b = i\circ j\in I\circ J$, it follows that
  \begin{align*}
       b + a 
       &= a^0 + b +a
       = -a + i\circ j + a
       = -a + i + \lambda_i\left(j\right) + a\\
       &= \left(-a + i + a\right) 
       + \left(-a + \lambda_i\left(j\right) + a\right)\in I\circ J.
   \end{align*}
   Thus, $I\circ J$ is a normal subset of $(S,+)$. Finally, we prove that $I\circ J$ is a normal subset of $(S,\circ)$. Now, observe that 
   $a^-= j^-\circ i^- = i^-\circ \left(i\circ j^-\circ i^-\right)\in I\circ J$ and, if $a, a^0\circ b\in I\circ J$, then
   \begin{align*}
      a\circ b = a\circ a^0\circ b
      = a + \lambda_a\left(a^0\circ b\right)\in I\circ J 
   \end{align*}
   and, if $a\circ b = i\circ j\in I\circ J$, 
  \begin{align*}
     b \circ a
     = a^- \circ (i\circ j )\circ a
     = \left(a^-\circ i\circ a\right)\circ 
     \left(a^-\circ j\circ a\right)
     \in I\circ J ,
  \end{align*}
  which completes the proof.
\end{proof}
\end{prop}

\medskip

In light of the results obtained in this section, it is interesting a further deepening of the dual weak brace theory. Our intent is to continue such an investigation in future work.
\medskip

 \section*{Acknowledgement}
We would like to thank the referee for the question related to \cref{+circ} that allows us to extend the result to weak braces that are not necessarily dual.

\bigskip

\bibliographystyle{elsart-num-sort}  
\bibliography{bibliography}

 \end{document}